\documentclass[12pt]{article}
\usepackage{epsfig,latexsym,amsfonts,amsmath,amssymb}
\newtheorem{theorem}{\textsf{Theorem}}[section]

\newtheorem{prop}{\textsf{Proposition}}[section]
\newtheorem{rem}{\textsf{Remark}}

\newcommand{\pro}{{\sc Proof}.~}
\newcommand{\R}{\mathbb{R}}
\def\Xint#1{\mathchoice
    {\XXint\displaystyle\textstyle{#1}}%
     {\XXint\textstyle\scriptstyle{#1}}%
     {\XXint\scriptstyle\scriptscriptstyle{#1}}%
     {\XXint\scriptstyle\scriptscriptstyle{#1}}%
	\!\int}
\def\XXint#1#2#3{{\setbox0=\hbox{$#1{#2#3}{\int}$}
	\vcenter{\hbox{$#2#3$}}\kern-.5\wd0}}

\begin{document}

\title{\normalsize 
\bf ON THE DIRICHLET PROBLEM \\
GENERATED BY THE MAZ'YA--SOBOLEV INEQUALITY}
\author { {\it A.I.~Nazarov}\footnote {Supported by grant NSh.4210.2010.1.},\\ 
Saint-Petersburg State University, \\
{\small e-mail:\ al.il.nazarov@gmail.com}} 
\date{}\maketitle

\section{Introduction}

In what follows, $x=(y;z)=(y_1,y';z)$ stands for a point in 
${\mathbb R}^n={\mathbb R}^m\times{\mathbb R}^{n-m}$, $n\ge3$, $2\le m\le n-1$.
Denote by ${\cal P}$ the subspace $\{x\in{\mathbb R}^n:\ y=0\}$; correspondingly,
${\cal P}^\bot=\{x\in{\mathbb R}^n:\ z=0\}$.

Let $\Omega$ be a domain in ${\mathbb R}^n$. By 
${\cal C}_0^{\infty}(\Omega)$ we denote the set of smooth functions with 
compact support in $\Omega$. For $1\le p<\infty$ we denote by $\dot W^1_p(\Omega)$ 
the closure of ${\cal C}_0^{\infty}(\Omega)$ with respect to 
the norm $\|\nabla v\|_{p,\Omega}$. Obviously, for bounded domains
$\dot W^1_p(\Omega)=\stackrel {o}{W}\!\vphantom{W}^1_p(\Omega)$.

By definition, for $0\le\sigma\le\min\{1,\frac np\}$ we put
$p^*_{\sigma}=\frac {np}{n-\sigma p}$.

\begin{prop}\label{prop}
The following inequality
\begin{equation}
\||y|^{\sigma-1}v\|_{p^*_{\sigma},\Omega}\le {\cal N}(p,\sigma,\Omega) 
\cdot \|\nabla v\|_{p,\Omega}.
\label{hardy}\end{equation}
holds true for any $v\in\dot W^1_p(\Omega)$ provided
\begin{equation}
\begin{array}{ll}
a)\ \Omega \mbox { is any domain in }\mathbb R^n&\mbox { for \ }\frac{n(p-m)}{p(n-m)}<\sigma\le1\quad
\mbox{(the region $\bf I$ on Fig. \ref{fig})};\\
b)\ \Omega \subset\mathbb R^n\setminus {\cal P}&\mbox { for \ }p>m,\ 
\sigma\le\min\{\frac{n(p-m)}{p(n-m)};\frac{n}{p}\},\ \sigma\ne1\ 
\mbox{($\bf II$ on Fig. \ref{fig})};\\
c)\ \Omega \subset\mathbb R^n\setminus (\ell\times\mathbb R^{n-m})& \mbox { for \ }p=m,\ \sigma=0\ 
\mbox{(black point on Fig. \ref{fig})}
\end{array}
\label{domain}\end{equation}
(here $\ell$ is a ray in ${\mathbb R}^m$ beginning at the origin).
\end{prop}


\begin{figure}\label{fig}
\begin{picture}(300,300)(-100,0)

\thinlines

\put(30,30){\vector(1,0){260}}
\put(30,30){\vector(0,1){260}}

\multiput(30,270)(5,0){12}{\line(1,0){2}}
\multiput(90,30)(0,5){48}{\line(0,1){2}}

\put(267,15){\large $1$}
\put(17,265){\large $1$}
\put(125,14){\large $\frac 1m$}
\put(85,14){\large $\frac 1n$}

\put(285,14){\large $\frac 1p$}
\put(17,285){\large $\sigma$}

\put(200,150){\Huge $\bf I$}
\put(67,90){\Huge $\bf II$}

\thicklines

\put(30,30){\line(1,0){240}}
\put(270,30){\line(0,1){240}}
\put(91,270){\line(1,0){179}}
\put(30,30){\line(1,4){59.6}}
\put(130,30){\line(-1,6){39.7}}
\put(90,270){\circle{3}} 
\put(130,30){\circle*{3}} 

\end{picture}
\caption{To the Proposition \ref{prop}}
\end{figure}


\pro
The case a) is well known; see, e.g., \cite[Sec.2.1.6]{M}. Note that for $\sigma=1$ we have classical Sobolev inequality.

Consider the cases b) and c). Note that it is sufficient to prove (\ref{hardy}) for 
$\Omega=\mathbb R^n\setminus {\cal P}$ (respectively, $\Omega=\mathbb R^n\setminus (\ell\times\mathbb R^{n-m})$).

For $\sigma=0$ one should take conventional Hardy inequality in 
$\mathbb R^m\setminus\{0\}$ (respectively, in $\mathbb R^m\setminus\ell$; see, e.g., \cite[Sec.2]{N})
and integrate it with respect to $z$. 

For $m<p<n$ the inequality (\ref{hardy}) can be obtained from the cases $\sigma=0$ and $\sigma=1$ 
by the H\"older inequality. For $p>n$ we also obtain (\ref{hardy}) by the H\"older inequality from the 
extreme cases $\sigma=0$ and $\sigma=\frac np$; the last one corresponds to the Morrey inequality, 
see \cite[Sec. 1.4.5]{M}.

Finally, we deal with the case $p=n$, $0<\sigma<1$. Consider the domain $\Omega=\Omega_1\times\Omega_2$, 
where $\Omega_1=B_2\setminus B_1\subset\mathbb R^m$ is a spherical layer while 
$\Omega_2=]0,1[^{n-m}\subset\mathbb R^{n-m}$ is a cube. Let us write down the embedding theorem
$W^1_n(\Omega)\hookrightarrow L_q(\Omega)$, with $q=\frac n{1-\sigma}$. Since the function $|y|^n$ 
is bounded and separated from zero in $\Omega$, this theorem can be rewritten as follows:
$$\Bigg(\int\limits_{\Omega}\frac {|v|^q}{|y|^n}\ dx\Bigg)^{n/q}\le 
C(q,m,n)\int\limits_{\Omega}\left(|\nabla v|^n+\frac {|v|^n}{|y|^n}\right)\ dx.$$
Note that all the terms in this inequality are invariant under translations in $z$ 
and under dilations in $x$. Therefore, the same inequality is valid for 
$\Omega_{k\mathfrak k}=2^k\big(\Omega_1\times(\Omega_2+\mathfrak k)\big)$, with 
$k\in\mathbb Z$, $\mathfrak k\in\mathbb Z^{n-m}$. Summing these inequalities we obtain, 
subject to $q>n$, 
$$\Bigg(\int\limits_{\mathbb R^n}\frac {|v|^q}{|y|^n}\ dx\Bigg)^{n/q}\le
\sum\limits_{k,\mathfrak k}\Bigg(\int\limits_{\Omega_{k\mathfrak k}}
\frac {|v|^q}{|y|^n}\ dx\Bigg)^{n/q}
\le C(q,m,n)\int\limits_{\mathbb R^n}\left(|\nabla v|^n+
\frac {|v|^n}{|y|^n}\right)dx.$$
The last term is already estimated, and we arrive at (\ref{hardy}).\hfill$\square$

\begin{rem}\label{rem} The assumption on $\Omega$ in the case c) can be considerably weakened. 
However, it is sharp for $\Omega$ being a wedge.
\end{rem}

We call (\ref{hardy}) {\bf the Maz'ya--Sobolev inequality}.\medskip

We are interested in the attainability of the sharp constant in (\ref{hardy}), i.e. in the attainability
of the norm of corresponding embedding operator. If $\Omega$ is unbounded, or 
$\overline\Omega\cap{\cal P}\ne\emptyset$, then this operator is, in general, noncompact; 
for $p<n$ and $\sigma=1$ this is the case for any $\Omega$. Therefore, the problem of attainability is 
nontrivial.

It is well known that the sharp constant in the Sobolev inequality ($p<n$ and $\sigma=1$)
does not depend on $\Omega$ and {\it is not attained} for any $\Omega$ provided the complement
of $\Omega$ is not negligible, i.e. $\dot W^1_p(\Omega)\ne\dot W^1_p(\mathbb R^n)$. We claim 
that the same is true for $p<n$ and $0<\sigma<1$ provided $\Omega\cap{\cal P}\ne\emptyset$. 
Indeed, since the inequality (\ref{hardy}) is dilation invariant, the sharp constant in this 
case cannot depend on $\Omega$ and equals ${\cal N}(p,\sigma,\mathbb R^n)$. Further, if 
the extremal function in (\ref{hardy}) exists, by standard argument (see, for example, the end 
of the proof of Theorem \ref{Th2}) it is (after a suitable normalization) a positive
generalized solution of the Dirichlet problem 
\begin{equation}\label{kkkk}
-\Delta_pu=\frac {u^{p^*_{\sigma}-1}}{|y|^{(1-\sigma)p^*_{\sigma}}} \quad
\mbox {in}\ \ \Omega;\qquad u\bigr|_{\partial\Omega}=0
\end{equation}
(here $\Delta_pu={\rm div}(|\nabla u|^{p-2}\nabla u)$ is $p$-Laplacian).

Extending $u$ by zero to $\mathbb R^n$, we obtain an extremal for (\ref{hardy}) in the whole space. Therefore,
this extension solves the equation (\ref{kkkk}) in $\mathbb R^n$, and thus it is positive in $\mathbb R^n$, a 
contradiction.\medskip 

By the way, it is worth to note that for $p=n$ the exponent in the denominator of 
(\ref{kkkk}) does not depend on $\sigma$ and equals $n$.\medskip

The case $\Omega\cap{\cal P}=\emptyset$, $\partial\Omega\cap{\cal P}\ne\emptyset$ is considerably 
more complicated. In the recent paper \cite{GhR1} the attainability of the sharp constant in 
(\ref{hardy}) was proved for $p=2$, $0<\sigma<1$, under rather restrictive assumptions on
(a smooth bounded domain) $\Omega$. Namely, it is supposed in \cite[Theorem 1.1]{GhR1} that all
the main curvatures at any point $x^0\in\partial\Omega\cap{\cal P}$ are nonpositive,
and the mean curvature at any such point does not vanish.\medskip

Our paper consists of two parts. First, we analyze the attainability of the sharp constant in 
(\ref{hardy}) for $\Omega$ being a wedge ${\cal K}=K\times\mathbb R^{n-m}$ (here $K$ is an open 
cone in $\mathbb R^m$) or a ``perturbed'' wedge. Here we consider all $1<p<\infty$ and 
$0\le\sigma<\min\{1,\frac np\}$. Naturally, we suppose that $\Omega$ satisfies (\ref{domain}).

In the second part we prove the attainability of the sharp constant in 
(\ref{hardy}) in a bounded domain for $p=2$ and $0<\sigma<1$ under considerably weakened 
requirements on $\partial\Omega$, see Section \ref{main} below. Unfortunately, we cannot transfer
this result to the case of arbitrary $p$ because we do not have in hands good estimates of
solutions to the model problem (\ref{kkkk}) in a half-space.\medskip

Let us discuss briefly the cases $m=1$ and $m=n$. For $m=1$ our problem of interest degenerates
in a sence\footnote{Note that Proposition \ref{prop} holds true for $m=1$ with the only exception:
the case c) should be attached to the case b). The proof runs without changes.}. Indeed, the 
only admissible wedge in this case is a half-space $\mathbb R^n_+=\{x\in\mathbb R^n:\ y_1>0\}$.
Theorems \ref{Th1} and \ref{Th2} in this case remain valid with the same proof while Theorems 
\ref{Th3} and \ref{Th4} are irrelevant. As for other domains, if $\Omega\subset\mathbb R^n_+$, and
$\partial\Omega\in{\cal C}^1$ touches $\cal P$, then in the neighborhood of a touching point 
$x^0$ $\Omega$ in the large scale looks like a half-space. Since (\ref{hardy}) is dilation 
invariant, we obtain ${\cal N}(p,\sigma,\Omega)\le {\cal N}(p,\sigma,\mathbb R^n_+)$. The reverse 
inequality is trivial. As in the case $\Omega\cap{\cal P}\ne\emptyset$, this implies 
non-attainability of the sharp constant in (\ref{hardy}) for any $\Omega$ provided the complement 
of $\Omega$ is not negligible in $\mathbb R^n_+$. For $p=2$ and bounded domain this fact was 
proved in \cite{GhR1}. Attainability of the sharp constant for $m=1$, $p=2$ in some unbounded 
domains without touching of $\cal P$ was discussed in \cite{Ti}\footnote{Example 1 after Lemma 2.7
in \cite{Ti} is not completely correct; it should be $\varphi>0$ instead of 
$\varphi\ge0$.}.\medskip

On the another hand, the problem for $m=n$, corresponding to the Hardy--Sobolev inequality, was
investigated in a number of papers. The existence of the extremal function in a cone was proved in 
\cite{N} (in the case $p=2$, $n\ge3$ this result was established earlier in \cite{E}). The problem 
in ``perturbed'' cone was considered in \cite{N1} (the case $p=2$, $\sigma=0$ was dealt with in 
\cite{PT}). For $\Omega$ being a compact Riemannian manifold with boundary, the conditions of 
attainability of the sharp constants in (\ref{hardy}) and in some similar inequalities were considered
in \cite{DN}. The case of bounded domains with $0\in\partial\Omega$ was treated in \cite{DN1} for $p=2$,
$n\ge2$; similar results under more restrictive assumptions on $\partial\Omega$ were obtained earlier 
in rather involved papers \cite{GhR} for $n\ge4$ and \cite{GhR2} for $n=3$. See also the survey \cite{N1},
where the history of related problems and extensive bibliography was given.\medskip

The paper is organized as follows. In Section \ref{wedges} we collect the results on existence
and qualitative properties of extremal functions in (\ref{hardy}) in wedges and in
wedges with compact perturbation bounded away from $\cal P$.

In Section \ref{main} we formulate the assumptions on the behavior of $\partial\Omega$ in 
a neighborhood of the origin and prove existence theorems for bounded domains. The technical 
estimates used in this proof are given in Sections \ref{denominator}--\ref{limit case}.\medskip

Let us introduce the following notation. 
${\mathbb S}_r^{n-1}$ is the sphere in ${\mathbb R}^n$ with radius $r$ centered at the origin; 
$\omega_{n-1}=\frac {2\pi^{n/2}}{\Gamma(\frac n2)}$ is the area of
${\mathbb S}_1^{n-1}$. 

We write $o_{\varepsilon}(1)$ to show the quantity tending to zero, as 
$\varepsilon\to0$, with other parameters assumed to be fixed. All the other 
$o(1)$ have the same meaning but are uniform with respect to $\varepsilon$.

We recall that a function $f: ]0,\delta[\to\mathbb R$ is 
{\bf regularly varying (RVF) of order $\alpha$} at the origin, 
if it has a constant sign, and for any $t>0$
$$\lim\limits_{\varepsilon\to0}\frac{f(\varepsilon t)}{f(\varepsilon)}=t^{\alpha}.
$$
For basic properties of RVFs see \cite{Se}.

We use letter $C$ to denote various positive constants. To indicate that $C$
depends on some parameters, we write $C(\dots)$.

\section{The Maz'ya--Sobolev inequality in wedges and in ``perturbed'' wedges}\label{wedges}

Our first statement provides the sharp constants in the Maz'ya inequality in wedges.

\begin{theorem}\label{Th1}
Let $2\le m\le n-1$, $1<p<\infty$, $\sigma=0$. Let $K$ be a cone in $\mathbb R^m$.
If $p\ge m$ we suppose that $K\ne {\mathbb R}^m$, and for $p=m$, in addition, 
$K\ne {\mathbb R}^m\setminus\{0\}$. Put $\Omega={\cal K}=K\times\mathbb R^{n-m}$ and
$G=K\cap{\mathbb S}_1^{m-1}$. Then the sharp constant in (\ref{hardy}) is not attained 
and equals $(\Lambda^{(p)}(G))^{-\frac 1p}$, where
\begin{equation}\label{p}
\Lambda^{(p)}(G)=\min_{v\in \stackrel {o}{W}\!\vphantom{W}^1_p(G)\setminus 
\{0\}}\frac {{\displaystyle\int_G}
\,\Bigl(\Bigl(\frac {m-p}p\Bigr)^2 v^2+|\nabla' v|^2\Bigr)^{\frac p2}dS}
{{\displaystyle\int_G}\,|v|^p\,dS}
\end{equation}
(here $\nabla'$ stands for the tangential gradient on 
${\mathbb S}_1^{m-1}\subset\mathbb R^m$).
\end{theorem}

\pro First, the minimum in (\ref{p}) is attained due to the compactness of 
embedding $\stackrel {o}{W}\!\!\vphantom{W}^1_p(G)\hookrightarrow L_p(G)$. 
Denote by $\widehat V$ the minimizer of (\ref{p}) normalized in $L_p(G)$. By standard 
argument, $\widehat V$ is positive in $G$.

Let us define $U(y,z)=U(y)=|y|^{1-\frac mp}\cdot\widehat V\big(\frac{y}{|y|}\big)$. 
It is shown in \cite[Theorem 18]{N1} that $U$ is a positive weak solution of the equation
\begin{equation}\label{pp}
-\Delta_p^{(y)}U=\Lambda^{(p)}(G)\ \frac {U^{p-1}}{|y|^p}\quad\mbox{in}\ K,
\qquad\mbox{and thus,}\qquad
-\Delta_pU=\Lambda^{(p)}(G)\ \frac {U^{p-1}}{|y|^p}\quad\mbox{in}\ {\cal K}.
\end{equation}

The relation $\Lambda^{(p)}(G)\le{\cal N}^{-p}(p,0,\Omega)$ follows now from 
\cite[Theorem 2.3]{PT1}. For the reader's convenience we reproduce the proof 
based on the so-called generalized Picone identity.

For any $u\in{\mathcal C}^{\infty}_0(\Omega)$ we set 
$h=\frac {|u|^p}{U^{p-1}}\in{\mathcal C}^1_0(\Omega)$. Then (\ref{pp}) implies
\begin{equation}\label{qq}
\gathered
\Lambda^{(p)}(G)\int\limits_{\Omega}\frac {|u|^p}{|y|^p}\,dx=
\Lambda^{(p)}(G)\int\limits_{\Omega}\frac {U^{p-1}}{|y|^p}h\,dx=
\int\limits_{\Omega}|\nabla U|^{p-2}\nabla U\cdot\nabla h\,dx=\\
=\int\limits_{\Omega}\Bigl(p|\nabla U|^{p-2}\nabla U\cdot\nabla u
\frac {|u|^{p-2}u}{U^{p-1}}-(p-1)|\nabla U|^p\frac {|u|^p}{U^p}\Bigr)\,dx
\stackrel {*}{\le}\\
\stackrel {*}{\le}\int\limits_{\Omega}\Bigl(p|\nabla u|\cdot|\nabla U|^{p-1}
\frac {|u|^{p-1}}{U^{p-1}}-(p-1)|\nabla U|^p\frac {|u|^p}{U^p}\Bigr)\,dx
\le\int\limits_{\Omega}|\nabla u|^p\,dx.
\endgathered
\end{equation}
Here $(*)$ is the Cauchy inequality while the last inequality follows from 
\begin{equation}\label{pppp}
r^p-prt^{p-1}+(p-1)t^p\ge0, \qquad r,t>0.
\end{equation}
By approximation, (\ref{qq}) holds true for 
$u\in\dot W^1_p(\Omega)$.\medskip

To prove $\Lambda^{(p)}(G)={\cal N}^{-p}(p,0,\Omega)$ we consider the sequence $u_\delta(y,z)=U_\delta(y)Z_\delta(z)$, where
$$U_\delta(y)=\begin{cases}
|y|^{1-\frac mp+\delta}\cdot\widehat V\big(\frac{y}{|y|}\big),& |y|\le R,\\
R^{1-\frac mp+\delta}\big(2-\frac {|y|}R\big)\cdot\widehat V\big(\frac{y}{|y|}\big),& R\le |y|\le 2R,\\
0,& |y|\ge 2R;\end{cases}\qquad
Z_\delta(Z)=\begin{cases}1,& |z|\le R,\\
2-\frac {|z|}{R},& R\le |z|\le 2R,\\
0,& |z|\ge 2R.\end{cases}
$$
Clearly, $u_\delta\in\dot W^1_p(\Omega)$. Direct 
computation shows
$$\int\limits_{\Omega}|\nabla u_\delta|^p\,dx=
\int\limits_{\Omega}\frac {|u_\delta|^p}{r^p}\,dx\cdot\big(\Lambda^{(p)}(G)+O(\delta)\big),
$$
and the statement follows.

Finally, the equality sign in $(*)$ means $\nabla u\parallel\nabla U$ while 
the equality in (\ref{pppp}) means $r=t$. These two facts imply
$$\frac {\nabla u}u=\frac {\nabla U}U\qquad \Longrightarrow \qquad u=cU$$ 
on the set $\{u\ne0\}$ and, therefore, in the whole $\Omega$. Since 
$U\!\notin\,\dot W^1_p(\Omega)$, the equality in (\ref{qq}) is 
impossible.\hfill$\square$\medskip


Next, we consider the Maz'ya--Sobolev inequality in wedges.

\begin{theorem}\label{Th2} Let $2\le m\le n-1$, $1<p<\infty$, 
$0<\sigma<\min\{1,\frac np\}$. Let $K$ be a cone in $\mathbb R^m$. If $p>m$ and
$\sigma\le\frac{n(p-m)}{p(n-m)}$ we suppose in addition that $K\ne {\mathbb R}^m$.
Put $\Omega={\cal K}=K\times\mathbb R^{n-m}$. Then the sharp constant in (\ref{hardy})
is attained, i.e. there exists a function $V\in\dot W^1_p(\Omega)$, $V>0$ in $\Omega$,
such that the inequality (\ref{hardy}) becomes equality.
\end{theorem}

\pro It is evident that the sharp constant in (\ref{hardy}) satisfies the relation
\begin{equation}
{\cal N}^{-1}(p,\sigma,\Omega)=\inf_{v\in \dot W^1_p({\Omega})\setminus\{0\}} J(v)
\equiv\inf_{v\in \dot W^1_p({\Omega})\setminus\{0\}}
\frac {\|\nabla v\|_{p,\Omega}}{\||y|^{\sigma-1}v\|_{p^*_{\sigma},\Omega}}.
\label{*}\end{equation}

Let $\{v_k\}$ be a minimizing sequence for the functional $J$.
Without loss of generality we can assume $\||y|^{\sigma-1}v_k\|_{p^*_{\sigma},\Omega}=1$
and $v_k\rightharpoondown v$ in $\dot W^1_p(\Omega)$. By the concentration-compactness 
principle of Lions (\cite{Ls}; see also \cite[Ch.1]{Ev}) we have
\begin{multline*}
||y|^{\sigma-1}v_k|^{p^*_{\sigma}}\rightharpoondown ||y|^{\sigma-1}v|^{p^*_{\sigma}}+
\sum\limits_{j\in{\cal M}}\alpha_j\delta(x-x^j), \\
|\nabla v_k|^p\rightharpoondown {\mathfrak M}\ge |\nabla v|^p+ {\cal N}^{-p}(p,\sigma,\Omega)
\sum\limits_{j\in{\cal M}}\alpha_j^{p/p^*_{\sigma}}\delta(x-x^j),
\end{multline*}
where the convergence is understood in the sense of measures on the one-point compactification
$\overline\Omega\cup\{\infty\}$, a set $\cal M$ is at most countable and $\alpha_j>0$.
Moreover, since the embedding $\dot W^1_p(\Omega)\hookrightarrow L_{p^*_{\sigma}}(\Omega)$ is
locally compact, we conclude that $x^j\in{\cal P}\cup\{\infty\}$.

Since $\{v_k\}$ is a minimizing sequence, by verbatim repetition of arguments 
from Theorem 2.2 \cite{LPT} we obtain the alternative~--- either $v_k\to v$ in 
$\dot W^1_p(\Omega)$ and $\cal M=\emptyset$ (in this case $v$ is a minimizer of $J$),
or $v=0$, $\cal M$ is a singleton and $\alpha=1$.

Let us remark here that, by the dilation invariance of the functional $J$, we can ensure the
additional relation 
$\int\limits_{\Omega\cap B_1}||y|^{\sigma-1}v_k|^{p^*_{\sigma}}\,dx=\frac12$, which
takes away the second variant.

It remains to note that the function $V=|v|$ also provides the minimum in the problem
(\ref{*}). Thus, after multiplying by a suitable constant, $V$ becomes a nonnegative
generalized solution of the Dirichlet problem to the Euler--Lagrange equation (\ref{kkkk})
and thus, it is super-$p$-harmonic in $\Omega$. By the Harnack inequality for $p$-harmonic functions
(see, e.g., \cite{Tr}), it is positive in $\Omega$.\hfill$\square$\medskip

Now we present some symmetry properties of the extremal function.

\begin{theorem}\label{Th3} Let the assumptions of Theorem \ref{Th2} be fulfilled. 
Then the functon $V$ providing the sharp constant in (\ref{hardy}) has the following properties:
\begin{enumerate}
 \item $V$ is radially symmetric with respect to $z$, i.e. $V=V(y;|z|)$;

 \item If $K$ is a circular cone, then $V$ is radially symmetric with respect to $y'$ and $z$, 
i.e. $V=V(y_1,|y'|;|z|)$;

 \item If $K=\mathbb R^m$ and $\sigma>\frac{n(p-m)}{p(n-m)}$, then $V$ is radially symmetric 
with respect to $y$ and $z$, i.e. $V=V(|y|;|z|)$;

 \item Let $K=\mathbb R^m\setminus\{0\}$. There exists $\widehat p\in\,]m,n[$, and for 
$p>\widehat p$ the function $\widehat\sigma(m,n,p)$ is defined, such that 
$\widehat\sigma<\min\{1,\frac np\}$ and for $\sigma>\widehat\sigma$ the function $V$ 
is not radially symmetric w.r.t. $y$.

\end{enumerate}

\end{theorem}

\pro 
1. This statement follows from the properties of the Schwarz symmetrization 
with respect to $z$-variables (or from the properties of the 
Steiner symmetrization with respect to $z_1$ for $m=n-1$). Indeed, this 
transformation does not enlarge the numerator in (\ref{*}), see, e.g.,
\cite[Ch.7]{PS}, and evidently retains the denominator.
Thus, it is sufficient to take infimum in (\ref{*}) over the set of functions,
radially symmetric w.r.t. $z$. Further, by the Euler equation (\ref{kkkk})
all critical points of an extremal radially symmetric w.r.t. $z$ have to be 
located at $\cal P$. In this case the numerator in (\ref{*}) strictly
decreases under symmetrization (see \cite{BZ}), and therefore no function 
asymmetric w.r.t. $z$ can provide the minimum in (\ref{*}).\medskip

2. In addition to the Part 1, in this case we can apply spherical symmetrization
along the spheres ${\mathbb S}_r^{m-1}$, which does not enlarge the numerator, 
see, e.g., \cite[App.C]{PS}, and retains the denominator.\medskip

3. Here we can apply the Schwarz symmetrization with respect to $y$-variables
 which does not enlarge the numerator, and does not reduce the denominator, see, e.g., 
\cite[Ch.3]{LL}.\medskip

4. In this case the Schwarz symmetrization in $y$s does not work, and we show that
the minimizer in general does not inherit the symmetry of extremal problem.\medskip

Let $u(|y|;|z|)$ be a function providing the minimum to the functional 
$J$ over the set of functions in $\dot W^1_p(\Omega)$, radially symmetric w.r.t. $y$ and $z$.
Without loss of generality, we assume that $\||y|^{\sigma-1}u\|_{p^*_{\sigma},\Omega}=1$. 
By the principle of symmetric criticality, see \cite{P}, $dJ_\sigma(u;h)=0$ 
for any variation $h\in\dot W^1_p(\Omega)$. 

Similarly to \cite[Theorem 1.3]{N2}, the second differential of $J$ at 
the point $u$ can be written as follows:
\begin{equation}\label{2deriv}
\gathered
J^{p-1}(u)\cdot d^2J(u;\ h) 
=\int\limits_{\Omega}|\nabla u| ^{p-4}\bigl((p-2)\langle
\nabla u, \nabla h \rangle^2 +|\nabla u|^2| \nabla h |^2\bigr) \,dx -\\
-J^p(u)\cdot\Bigl[(p-p^*_{\sigma}) \cdot \Bigl(\,\int\limits_{\Omega}
\frac {|u|^{p^*_{\sigma}-2}u h}{|y|^{(1-\sigma)p^*_{\sigma}}}\,dx
\Bigr)^2+(p^*_{\sigma}-1)\cdot\!\int\limits_{\Omega}\frac 
{|u|^{p^*_{\sigma}-2} h^2}{|y|^{(1-\sigma)p^*_{\sigma}}}\,dx\Bigr].
\endgathered
\end{equation}

Now we set $h(y;z)=u(|y|;|z|)\cdot\frac{y_1}{|y|}$. By symmetry of 
$u$, $\int\limits_{\Omega}\frac {|u|^{p^*_{\sigma}-2}u h}
{|y|^{(1-\sigma)p^*_{\sigma}}}\, dx=0$. Substituting into 
(\ref{2deriv}), we obtain
$$  J^{p-1}(u)\cdot d^2J(u;\ h)
=\int\limits_{\Omega} |\nabla u|^{p-2} \frac {u^2}{|y|^2}\,dx
-J^p(u)\cdot \frac{p^{*}_\sigma-p}{m-1}\cdot \int\limits_{\Omega} \frac
{|u|^{p^{*}_\sigma} f^2}{|y|^{(1-\sigma)p^{*}_{\sigma}}}\,dx.
$$
Finally, we estimate the first integral by H\"older and Hardy 
inequalities and arrive at
$$d^2J(u;\ h)\le J(u)\cdot\left[\left(
\frac p{p-m}
\right)^2-
\frac {p^2\sigma}{(m-1)(n-p\sigma)}
\right].$$
If $p\ge n$ then the quantity in square brackets is negative for $\sigma$ close to 
$\frac np$. If $p<n$ is close to $n$, this quantity is also negative for $\sigma$ 
close to $1$. In both cases the statement follows.\hfill$\square$\medskip

{\bf Corollary}. For $p>\widehat p$ and $\widehat\sigma<\sigma<\min\{1,\frac np\}$ 
the problem (\ref{kkkk}) in $\mathbb R^n\setminus{\cal P}$ has at least two 
nonequivalent positive solutions.\medskip

\pro The first solution is a global minimizer of $J$ (under suitable normalization),
the second one is a minimizer over the set of functions symmetric w.r.t. $y$.\hfill$\square$\medskip

Further, we consider $\Omega$ being a perturbed wedge.

\begin{theorem}\label{Th4}
Suppose that $2\le m\le n-1$, $1<p<\infty$ and $0\le\sigma<\min\{1,\frac np\}$. 
Let $\Omega_1={\cal K}=K\times\mathbb R^{n-m}$ be a wedge satisfying (\ref{domain}),
$\Omega_2\Subset{\mathbb R}^n\setminus{\cal P}$ and $\Omega_1\cap\Omega_2\ne\emptyset$.

{\bf 1}. For $\Omega=\Omega_1\setminus\overline\Omega_2$ 
is not attained.

{\bf 2}. Let $\sigma>0$. Then for $\Omega=\Omega_1\cup\Omega_2$ the sharp constant in (\ref{hardy}) is
attained provided $\dot W^1_p(\Omega)\ne\dot W^1_p(\Omega_1)$.
 
{\bf 3}. Let $\sigma=0$. Then, given $\Omega'_2\Subset{\mathbb R}^m\setminus\{0\}$, $\Omega'_2\cap K\ne\emptyset$,
there exists $L<\infty$ such that if $\Omega_2\supset\Omega'_2\times\,]-L,L[$, $\Omega=\Omega_1\cup\Omega_2$ 
and $\dot W^1_p(\Omega)\ne\dot W^1_p(\Omega_1)$ then the sharp constant in (\ref{hardy}) is attained.
\end{theorem}

\pro 
{\bf 1}. For any $u\in{\mathcal C}_0^\infty({\Omega_1})$ there exists a 
dilation $\Pi$ such that $\Pi u\in{\mathcal C}_0^\infty({\Omega})$. Due to 
the dilation invariance of (\ref{hardy}) we conclude that
${\cal N}(p,\sigma,\Omega)={\cal N}(p,\sigma,\Omega_1)$.

Thus, if $u$ minimizes the quotient (\ref{*}) on $\dot W^1_p(\Omega)$ then 
its zero continuation minimizes (\ref{*}) on $\dot W^1_p(\Omega_1)$. 
Therefore, it is the nonnegative solution of the problem (\ref{kkkk}) in $\Omega_1$. By 
Harnack's inequality for $p$-harmonic functions, it is positive in $\Omega_1$, a 
contradiction.\medskip

{\bf 2}. By Theorem \ref{Th2}, there exists a function $u$ positive in $\Omega_1$ 
that minimizes the quotient (\ref{*}) on $\dot W^1_p(\Omega_1)$. If 
${\cal N}(p,\sigma,\Omega)={\cal N}(p,\sigma,\Omega_1)$ then the zero 
continuation of $u$ minimizes (\ref{*}) on $\dot W^1_p(\Omega)$ that again 
leads to contradiction. Therefore, ${\cal N}(p,\sigma,\Omega)>{\cal N}(p,\sigma,\Omega_1)$. 

Now the statement follows by the concentration-compactness principle. Indeed, let $\{v_k\}$ 
be a minimizing sequence for the functional $J$. Without loss of generality we can assume
$\||y|^{\sigma-1}v_k\|_{p^*_{\sigma},\Omega}=1$ and $v_k\rightharpoondown v$ in $\dot W^1_p(\Omega)$. 
As in Theorem \ref{Th2}, if $v_k\not\to v$ then
$$||y|^{\sigma-1}v_k|^{p^*_{\sigma}}\rightharpoondown \delta(x-\widehat x), \qquad
|\nabla v_k|^p\rightharpoondown {\cal N}^{-p}(p,\sigma,\Omega)\,\delta(\widehat x-x),
$$
and $\widehat x\in{\cal P}\cup\{\infty\}$. 

Since $\Omega_2\Subset{\mathbb R}^n\setminus{\cal P}$, similarly to the proof of Corollary 2.1 \cite{LPT}, 
we can assume that $v_k\big|_{\Omega_2}\equiv0$. This implies 
${\cal N}(p,\sigma,\Omega)\le{\cal N}(p,\sigma,\Omega_1)$, a contradiction.\medskip 

{\bf 3}. Define $\widetilde\Omega'=K\cup\Omega'_2$ and $\widetilde\Omega=\widetilde\Omega'\times\mathbb R^{n-m}$.
It is proved in \cite[Theorem 20]{N1} that ${\cal N}(p,0,\widetilde\Omega')>{\cal N}(p,0,K)$, and there 
exists a minimizer $\widetilde U$ of the quotient (\ref{*}) in $\widetilde\Omega'$. Then
$\widetilde U$ is a positive weak solution of the equation
\begin{equation*}
-\Delta_p^{(y)}\widetilde U={\cal N}^{-p}(p,0,\widetilde\Omega')\ \frac {\widetilde U^{p-1}}{|y|^p}\quad\mbox{in}\ \widetilde\Omega',\qquad\mbox{and thus,}\qquad
-\Delta_p\widetilde U={\cal N}^{-p}(p,0,\widetilde\Omega')\ \frac {\widetilde U^{p-1}}{|y|^p}\quad\mbox{in}\ \widetilde\Omega.
\end{equation*}
As in Theorem \ref{Th1}, this implies
$${\cal N}(p,0,\widetilde\Omega)\ge{\cal N}(p,0,\widetilde\Omega')>{\cal N}(p,0,K)={\cal N}(p,0,\Omega_1)$$
(the last equality is due to Theorem \ref{Th1}).

Thus, there exists $u\in{\mathcal C}_0^\infty(\widetilde\Omega)$ such that 
$\||y|^{-1}u\|_{p,\widetilde\Omega}>{\cal N}(p,0,\Omega_1) \cdot \|\nabla u\|_{p,\widetilde\Omega}$.
This means ${\cal N}(p,0,\Omega)>{\cal N}(p,0,\Omega_1)$ if $L$ is sufficiently large, and the statement
follows by the concentration-compactness principle.\hfill$\square$\medskip

In what follows we need some estimates for the solution of the extremal problem (\ref{*}) for $p=2$
in the half-space. For the sake of brevity, we denote 
$$q=2^*_{\sigma}=\frac{2n}{n-2\sigma};\qquad \mu_q(\Omega)={\cal N}^{-1}(2,\sigma,\Omega);\qquad
\mu_q=\mu_q(\mathbb R^n_+).
$$
By $\phi$ we denote a minimizer of the problem (\ref{*}) for $p=2$ in $\Omega={\mathbb R}^n_+$. 
Without loss of generality we can assume $\||y|^{\sigma-1}\phi\|_{q,{\mathbb R}^n_+}=1$.
Then $\phi$ is a weak solution of the Dirichlet problem
\begin{equation}
-\Delta u=\mu^2_q\cdot\frac{u^{q-1}}{|y|^{q(1-\sigma)}}\quad\mbox{in}
\quad{\mathbb R}^n_+, \qquad u\big|_{x_n=0}=0.
\label{eq}\end{equation}

\begin{prop}\label{prop1} The function $\phi$ satisfies the following relations:
\begin{equation}
\phi(x)\sim Cx_n,\qquad
|\nabla \phi(x)|\sim C,\qquad x\to0;
\label{estimate1}\end{equation}
\begin{equation}
\phi(x)\sim\frac{Cx_n}{|x|^{n}},\qquad
|\nabla \phi(x)|\asymp\frac{C}{|x|^n},\qquad x\to\infty.
\label{estimate2}\end{equation}
\end{prop}

\pro First, we claim that $\phi\in {\cal C}^{1+\gamma}_{\rm{loc}}(\overline{{\mathbb R}^n_+})$.
Indeed, the standard elliptic theory, see, e.g., \cite{LU}, provides 
$\phi\in {\cal C}^2_{\rm{loc}}(\overline{{\mathbb R}^n_+}\setminus{\cal P})$.
Estimates in the neighborhood of $\cal P$ can be obtained using elliptic theory in domains
with edges, see, e.g., \cite{MP}. Note that the property 
$\phi\in {\cal C}^1_{\rm{loc}}(\overline{{\mathbb R}^n_+})$ was proved also in 
\cite[Appendix]{GhR1}. 

Further, the Hopf lemma gives $\phi_{x_n}\big|_{x_n=0}>0$, and 
(\ref{estimate1}) follows.

Finally, the relations (\ref{estimate2}) follow from (\ref{estimate1}). 
Indeed, the direct computation shows that the image of 
$\phi$ under the Kelvin transform is also a solution of the problem 
(\ref{eq}) while (\ref{estimate1}) turns into (\ref{estimate2}).\hfill$\square$

\section{The case of bounded domain}\label{main}

We assume that in a neighborhood of the set ${\cal P}\cap\partial\Omega$
the boundary is of class ${\cal C}^1$; outside this neighborhood we impose no 
assumptions on $\partial\Omega$. Suppose there exists a point
$x^0\in{\cal P}\cap\partial\Omega$ (without loss of generality, $x^0=0$)
satisfying the properties listed below.

Let us introduce local Cartesian coordinates with $y'=(y_2,\dots,y_m)$ 
in the tangent plane and the axis $Oy_1$ directed into $\Omega$. Then in a 
neighborhood of the origin $\partial\Omega$ is given by equation $y_1=F(y';z)$. 
It is evident that $F\in {\cal C}^1$ and $F(y';z)=o(|y'|+|z|)$. Moreover, the assumption
${\cal P}\cap\Omega=\emptyset$ implies $F(0;z)\ge0$.

We say that $\partial\Omega$ is {\bf average concave} in a neighborhood of the origin
(see \cite{DN1}), if for sufficiently small $\rho$
\begin{equation}
f(\rho):=\!\!\Xint{\ \ -}\limits_{{\mathbb S}_\rho^{n-2}}\!\!F(y';z)\,d{\mathbb S}_\rho(y',z)<0
\label{eq1}\end{equation}
(here and later the dashed integral stands for the mean value).

We introduce also the functions
$$f_1(r;t):=\!\!\Xint{\ \ \,-}\limits_{{\mathbb S}_r^{m-2}}\!\!\!\Xint{\quad\ \ -}\limits_{{\mathbb S}_t^{n-m-1}}\!\!F(y';z)\,
d{\mathbb S}_r(y')d{\mathbb S}_t(z),$$
$$f_2(\rho):=\!\!\Xint{\ \ -}\limits_{{\mathbb S}_\rho^{n-2}}\!\!|\nabla'F(y';z)|^2\,d{\mathbb S}_\rho(y',z),$$
($\nabla'$ stands for the gradient with respect to $(y',z)$) and assume that for sufficiently small $\rho$
\begin{equation}
\cos^{m-2}(\beta)\sin^{n-m-1}(\beta)\cdot|f_1(\rho\cos(\beta),\rho\sin(\beta))|\le C\cdot|f(\rho)|,\qquad
\beta\in[0,\frac{\pi}{2}],
\label{eq11}\end{equation}
and
\begin{equation}
\lim\limits_{\rho\to 0}\frac{f_2(\rho)}{f(\rho)}\,\rho=0.
\label{eq4}\end{equation}

We say that $\partial\Omega$ is {\bf average concave in ${\cal P}$ and ${\cal P}^\bot$ directions}
in a neighborhood of the origin, if (\ref{eq1}) holds for sufficiently small $\rho$, and 
\begin{equation}
\Phi(\beta):=\lim\limits_{\rho\to0}\frac{f_1(\rho\cos(\beta),\rho\sin(\beta))}{f(\rho)}\ge0,\qquad
\beta\in[0,\frac{\pi}{2}]. 
\label{eq2}\end{equation}


Now we can formulate the main result of the second part of our paper.

\begin{theorem}\label{Th5} Let $\partial\Omega$ be average concave in ${\cal P}$ and 
${\cal P}^\bot$ directions in a neighborhood of the origin, and let the relations (\ref{eq11}) 
and (\ref{eq4}) hold. Suppose also that $f$ is regularly varying of order $\alpha\in[1, n+1[$
at the origin. Then for $p=2$ and for any $0<\sigma<1$ the infimum in (\ref{*}) is attained.
\end{theorem}

Let us compare our assumptions with those of \cite{GhR1}. If $\partial\Omega$ is smooth and $\alpha=2$, then
$$f(\rho)\sim {\cal H}\rho^2, \qquad f_1(r,t)\sim {\cal H}^{\cal P}r^2+{\cal H}^{{\cal P}^\bot}t^2,
\qquad f_2(\rho)\sim C\rho^2$$ 
near the origin (here ${\cal H}=\frac{1}{2(n-1)}{\bf Sp}(\nabla^{\prime2}F(0))$ is the mean curvature of 
$\partial\Omega$ at the origin; respectively, ${\cal H}^{\cal P}=\frac{1}{2(m-1)}{\bf Sp}(\nabla_{y'}^2F(0))$ 
and ${\cal H}^{{\cal P}^\bot}=\frac{1}{2(n-m)}{\bf Sp}(\nabla_z^2F(0))$.

Since ${\cal P}\cap\Omega=\emptyset$, ${\cal H}^{{\cal P}^\bot}$ is always non-negative. Thus, the relations 
(\ref{eq1}) and (\ref{eq2}) mean that 
\begin{equation}
{\cal H}^{\cal P}<0;\qquad {\cal H}^{{\cal P}^\bot}=0.
\label{meancurv}\end{equation}
The relations (\ref{eq11}) and (\ref{eq4}) are automatically fulfilled in this case. 

One can see that (\ref{meancurv}) is considerably weaker then the assumptions of \cite[Theorem 1.1]{GhR1}. 
We underline also that our hypotheses must be fulfilled at {\bf some} point $x^0\in{\cal P}\cap\partial\Omega$
while the authors of \cite{GhR1} constrain the curvatures at {\bf any} point $x^0\in{\cal P}\cap\partial\Omega$.
Moreover, we do not require even the existence of the mean curvature (if $\alpha<2$). On the 
other hand, for $\alpha>2$ all curvatures vanish at the origin. 

\begin{rem}\label{rem2}
The assumption (\ref{eq4}) can fulfil even if the main term of the 
asymptotic expansion of $F$ vanishes under average. For example, it is the 
case if $F(y';z)=y_2^3-y_3^4$.
\end{rem}

\begin{rem}\label{rem1}
The assumption (\ref{eq11}) is used only to ensure the limit passage under integral sign
and can be easily weakened. However, it cannot be removed at all, and we prefer to give
it in a simple form. In turn, the assumption (\ref{eq2}) could be weakened if we had in hands
more detailed information on the function $\phi$.
\end{rem}

Now consider the limit case $\alpha=n+1$. In this case we can drop the assumption (\ref{eq2}).

\begin{theorem}\label{Th6} Let $\partial\Omega$ be average concave in a neighborhood of the origin, 
and let the relations (\ref{eq11}) and (\ref{eq4}) hold. Suppose also that $f$ is regularly varying 
of order $n+1$ at the origin, and 
$\int_0^\delta\frac{f(\mathfrak r)}{\mathfrak r^{n+2}}\,d\mathfrak r=-\infty$.
Then for $p=2$ and for any $0<\sigma<1$ the infimum in (\ref{*}) is attained.
\end{theorem}

{\sc Proof of Theorems \ref{Th5} and \ref{Th6}}.~ Let $\{v_k\}$ be a minimizing sequence for (\ref{*}).
Without loss of generality we can assume $\||y|^{\sigma-1}v_k\|_{q,\Omega}=1$ and 
$v_k\rightharpoondown v$ in $\dot W^1_2(\Omega)$. 

Operating as in the proof of Theorem \ref{Th2}, we obtain the alternative~--- either $v$ is a
a minimizer of the extremal problem, or $v=0$ and
$$||y|^{\sigma-1}v_k|^q\rightharpoondown \delta(x-\widehat x), \qquad
|\nabla v_k|^2\rightharpoondown \mu_q^2(\Omega)\delta(x-\widehat x), \qquad
\widehat x\in{\cal P}\cap\partial\Omega$$
(the convergence is understood in the sense of measures on $\overline\Omega$).

We claim that in the second case $\mu_q(\Omega)\ge\mu_q$. Indeed, without loss of generality, 
$v_k$ concentrate near the origin. Further, as in the Corollary 2.1 \cite{LPT}, we can assume 
supports of $v_k$ located in arbitrarily small ball. Since $F(y';z)=o(|y'|+|z|)$ and $F(0;z)\ge0$, 
this implies
$${\rm supp}(v_k)\subset {\cal K}_{\varkappa}:=\{x\in{\mathbb R}^n:\
y_1>-\varkappa|y'|\}$$
for any $\varkappa>0$. Hence
$$\mu_q(\Omega)\ge\lim\limits_{\varkappa\to0}\mu_q({\cal K}_{\varkappa})=
\mu_q({\cal K}_0)=\mu_q
.$$

Therefore, to prove the statements we need only to produce a 
function having the quotient (\ref{*}) less then $\mu_q$. Similarly to 
\cite{DN1}, we construct such function using a suitable dilation and 
``bending'' of the function $\phi$ and multiplying it by a cut-off function 
with small support. The sharp estimates of behavior of $\phi$ (Proposition
\ref{prop1}) provide the desired result under assumptions on $\partial\Omega$ close to optimal.\medskip

Choose $\delta$ such that for $|y'|+|z|<2\,\delta$ the relation (\ref{eq1}) is 
satisfied and $|F(y';z)|\le \frac{|y'|+|z|}{2}$.

Let us introduce the coordinate transformation $\Theta_\varepsilon: x\mapsto
\varepsilon^{-1}(x-F(y';z)\,e_m)$. It is evident that in a neighborhood of the 
origin $\Theta_\varepsilon$ straightens $\partial\Omega$; its Jacobian 
equals $\varepsilon^{-n}$. Also it is easy to see that for $r<\delta$ we have
$B_{\frac{r}{2\,\varepsilon}}\subset \Theta_\varepsilon(B_r)\subset
B_{\frac{2\,r}{\varepsilon}}$.

Let $\widetilde{\varphi}\in {\cal C}_0^\infty(\R^n)$ be a function, radially symmetric
w.r.t. $y$ and $z$ and satisfying $0\le\widetilde{\varphi}\le1$,
$$\widetilde{\varphi}(x)\equiv\widetilde{\varphi}(|y|;|z|)=
\left\{%
\begin{array}{ll}
    1, & \mbox {if $|y|<\frac{\delta}{2}$ and $|z|<\frac{\delta}{2}$;} \\
    0, & \mbox {if $|y|>\delta$ or $|z|>\delta$.} \\
\end{array}%
\right.
$$
We introduce the cut-off function $\varphi(x)=\widetilde{\varphi}(\Theta_1(x))$.
Obviously, the function $x\mapsto \varphi(\Theta_\varepsilon^{-1}(x))$ is 
radially symmetric w.r.t. $y$ and $z$:
$$\varphi(\Theta_\varepsilon^{-1}(x))=\varphi(\varepsilon\,y_1+F(\varepsilon\,y';\varepsilon z)
,\varepsilon\,y';\varepsilon z)=\widetilde{\varphi}
(\varepsilon\,y_1,\varepsilon\,y'; \varepsilon z)=\widetilde{\varphi} (\varepsilon\,|y|; \varepsilon\,|z|).$$

Now we define the function
$$\phi_\varepsilon(x)=\varepsilon^{-(n-2)/2}
\phi\Big(\Theta_\varepsilon(x)\Big)\,\varphi(x).$$
It is easy to see that
$\phi_\varepsilon\in \stackrel {o}{W}\!\vphantom{W}^1_2(\Omega)$, 
if $\delta$ and $\varepsilon$ are sufficiently small.

In Sections \ref{denominator}--\ref{limit case} we show that
\begin{equation}
\int\limits_{\Omega}\frac{|\phi_\varepsilon(x)|^{q}}{|y|^{q(1-\sigma)}}\,dx=
1-A_1(\varepsilon)(1+o_\delta(1)+o_\varepsilon(1)),
\label{denom}\end{equation}
\begin{equation}
\int\limits_{\Omega}|\nabla\phi_\varepsilon(x)|^2\,dx=
\mu^2_q+A_2(\varepsilon)(1+o_\delta(1)+o_\varepsilon(1))-
\frac {2\mu^2_q}{q}A_1(\varepsilon)(1+o_\varepsilon(1))
\label{numer}\end{equation}
(we recall that $o_{\delta}(1)$ is uniform with respect to $\varepsilon$). 
For given $\delta$, in these formulas we have, as $\varepsilon\to0$,
\begin{equation}
A_1(\varepsilon)\sim C\varepsilon^{-1}\,f(\varepsilon);
\label{A1}\end{equation}
\begin{equation}
A_2(\varepsilon)\sim \left\{%
\begin{array}{ll}
    C\varepsilon^{-1}\,f(\varepsilon), & \hbox{under assumptions of Theorem \ref{Th5};} \\ \\
    C\varepsilon^{n}\int_\varepsilon^\delta\frac{f(\mathfrak r)}{\mathfrak r^{n+2}}\,d\mathfrak r, &
    \hbox {under assumptions of Theorem \ref{Th6}.} \\
\end{array}%
\right.
\label{A2}\end{equation}

The relations (\ref{A1}) and (\ref{A2}) imply 
$A_1(\varepsilon)=O(A_2(\varepsilon))$ (in the case $\alpha=n+1$ it follows 
from (\ref{svf})). Therefore, for sufficiently small $\delta$ and 
$\varepsilon$ we have, subject to (\ref{eq1}),
 \begin{multline}
\frac {\|\nabla \phi_\varepsilon\|_2^2}
{\||y|^{\sigma-1}\phi_\varepsilon\|^2_q}=
\frac {\mu^2_q+A_2(\varepsilon)(1+o_\delta(1)+o_\varepsilon(1))-
\frac {2\mu^2_q}{q}A_1(\varepsilon)(1+o_\varepsilon(1))}
{\Bigl(1-A_1(\varepsilon)(1+o_\delta(1)+o_\varepsilon(1))\Bigr)^{2/q}}=\\
=\mu^2_q+A_2(\varepsilon)(1+o_\delta(1)+o_\varepsilon(1))<\mu^2_q, 
\nonumber\end{multline}
and both Theorems follow.\hfill$\square$

\section{Estimate of the denominator and derivation of (\ref{denom})}\label{denominator}

We have, using the Taylor expansion,
\begin{multline*}
\int\limits_{\Omega}\frac{|\phi_\varepsilon(x)|^{q}}{|x|^{q(1-\sigma)}}dx=
\int\limits_{\R^n_+} \frac{|\phi(y;z)|^{q}}
{|y+\varepsilon^{-1}\,F(\varepsilon y';\varepsilon z)\,e_m|^{q(1-\sigma)}}\,
\varphi^q(\Theta_\varepsilon^{-1}(x))\,dydz=\\
=\int\limits_{\R^n_+}\frac{|\phi(y;z)|^{q}}{|y|^{q(1-\sigma)}}
\,\widetilde\varphi^q(\varepsilon y;\varepsilon z)\cdot 
\bigg(1-q(1-\sigma)\,\frac{F(\varepsilon y';\varepsilon z)\,y_1}{\varepsilon|y|^2}\bigg)\,dydz+\\
+O_\delta(1)\int\limits_{\R^n_+} \frac{|\phi(y;z)|^{q}F^2(\varepsilon y';\varepsilon z)}
{\varepsilon^2|y+\xi\varepsilon^{-1}\,F(\varepsilon y';\varepsilon z)\,e_m|^{q(1-\sigma)+2}}\,
\widetilde\varphi^q(\varepsilon y;\varepsilon z)\,dydz=:I_1-I_2+I_3
\end{multline*}
(here $\xi=\xi(y,z)\in\,]0,1[$).

\begin{enumerate}
    \item Since $\phi$ is normalized, $I_1\le1$. On the another hand, the first estimate in
(\ref{estimate2}) gives
$$1-I_1=
\int\limits_{\R^n_+}\frac{\phi^{q}(y;z)}{|y|^{q(1-\sigma)}}
\Big(1-\widetilde\varphi^q(\varepsilon y;\varepsilon z)\Big)\,dydz
\le C\int\limits_{\R^n_+\setminus B_{\frac{\delta}{2\varepsilon}}}
\frac{|y|^{-q(1-\sigma)}dydz}{(|y|^2+|z|^2)^{\frac{n-1}{2}}}
\le C\,\Big(\frac{\varepsilon}{\delta}\Big)^{\frac{qn}{2}}.
$$
     \item 
$$I_2=\frac {q(1-\sigma)}{\varepsilon}\int\limits_{\R^n_+}
\frac{\phi^{q}(y;z)y_1}{|y|^{q(1-\sigma)+2}}\,
\widetilde{\varphi}^q(\varepsilon y;\varepsilon z)\ F(\varepsilon y';\varepsilon z)\,dydz
=:A_1(\varepsilon).$$

\begin{prop}\label{prop2} Given $\delta$, the function $A_1(\varepsilon)$ 
satisfies (\ref{A1}), as $\varepsilon\to0$.
\end{prop}
\pro We claim that
\begin{multline}
\lim\limits_{\varepsilon\to
0}\varepsilon\frac{A_1(\varepsilon)}{f(\varepsilon)}={q(1-\sigma)}\,
\omega_{m-2}\omega_{n-m-1}\times\\
\times\int\limits_0^{\infty}\rho^{n-2+\alpha}\int\limits_0^{\frac{\pi}{2}}
\cos^{m-2}(\beta)\sin^{n-m-1}(\beta)\Phi(\beta)\int\limits_0^{\infty}
\frac{\phi^{q}(\rho\cos(\beta),s;\rho\sin(\beta))\,s\,ds}
{(\rho^2\cos^2(\beta)+s^2)^{\frac{{q(1-\sigma)}+2}{2}}}\,d\beta d\rho.
\label{eq6}\end{multline}

To prove this we apply the Lebesgue theorem. We have
\begin{multline*}
\frac{\varepsilon A_1(\varepsilon)}{q(1-\sigma)\ f(\varepsilon)}=
\frac{1}{f(\varepsilon)}\int\limits_{\R_+^n}
\widetilde{\varphi}^q(\varepsilon y;\varepsilon z)\frac{\phi^{q}(y;z)y_1}
{|y|^{q(1-\sigma)+2}}\ F(\varepsilon y';\varepsilon z)\,dydz=\\
=\frac{1}{f(\varepsilon)}
\int\limits_0^{\infty}\int\limits_0^{\infty}
\int\limits_0^{\infty}\widetilde{\varphi}^q(\varepsilon\sqrt{r^2+y_1^2};\varepsilon t)
\frac{\phi^{q}(y_1,r;t)y_1}{(r^2+y_1^2)^{\frac {q(1-\sigma)+2}{2}}}\times\\
\times\int\limits_{{\mathbb S}_r^{m-2}}\int\limits_{{\mathbb S}_t^{n-m-1}}\!\!F(\varepsilon y';\varepsilon z)\,
d{\mathbb S}_r^{m-2}(y')\,d{\mathbb S}_t^{n-m-1}(z)\,dy_1\,drdt=\\
=\omega_{m-2}\omega_{n-m-1}\int\limits_0^{\infty}\int\limits_0^{\infty}
r^{m-2}t^{n-m-1}\frac{f_1(\varepsilon r;\varepsilon t)}{f(\varepsilon)}\times\\
\times\int\limits_0^{\infty}
\widetilde{\varphi}^q(\varepsilon\sqrt{r^2+s^2};\varepsilon t)
\frac{\phi^{q}(r,s;t)\,s\,ds}{(r^2+s^2)^{\frac{{q(1-\sigma)}+2}{2}}}\,drdt.
\end{multline*}

We can apply the Monotone Convergence Theorem to the interior integral. Further, the regular behavior of $f$
implies
$$\lim\limits_{\varepsilon\to0}
\frac{f_1(\varepsilon \rho\cos(\beta);\varepsilon \rho\cos(\beta))}{f(\varepsilon)}=
\lim\limits_{\varepsilon\to0}
\frac{f_1(\varepsilon \rho\cos(\beta);\varepsilon \rho\cos(\beta))}{f(\varepsilon\rho)}\cdot
\lim\limits_{\varepsilon\to0}
\frac{f(\varepsilon \rho)}{f(\varepsilon)}=\Phi(\beta)\rho^\alpha.
$$
Therefore, the assumption on the pointwise convergence is satisfied. Now we 
produce a summable majorant. 

Due to the estimates (\ref{estimate1}) and (\ref{estimate2}), the interior integral is 
bounded from above by
$$C\chi^{(\varepsilon)}(r;t)\int\limits_0^\infty
\frac{
s^{q\sigma-1}\,ds}{1+(r^2+s^2+t^2)^{\frac{qn}2}}
\le \frac{C\chi^{(\varepsilon)}(r;t)}{1+(t^2+r^2)^{\frac{q(n-\sigma)}{2}}},$$
where $\chi^{(\varepsilon)}(r;t)=\chi_{[0,\frac{\delta}{\varepsilon}]}(r)
\cdot\chi_{[0,\frac{\delta}{\varepsilon}]}(t)$.

Now we pass to the polar coordinates. Using (\ref{eq11}), we estimate the integrand by
$$C\chi_{[0,\frac{2\delta}{\varepsilon}]}(\rho)
\cdot\frac{f(\varepsilon \rho)}{f(\varepsilon)}\cdot\frac{\rho^{n-2}}{1+\rho^{q(n-\sigma)}}.
$$
For $\gamma>0$ we have
$$\chi_{[0,\frac{2\delta}{\varepsilon}]}(\rho)
\cdot\frac{f(\varepsilon \rho)}{f(\varepsilon)}
\le
\chi_{[0,1]}(\rho)\,\rho^{\alpha-\gamma}\frac{f(\varepsilon\rho)\,
(\varepsilon\, \rho)^{-\alpha+\gamma}}{f(\varepsilon)\,\varepsilon^{-\alpha+\gamma}}\,+
\chi_{[1,\frac{2\delta}{\varepsilon}]}(\rho)\,\rho^{\alpha+\gamma}
\frac{f(\varepsilon\rho)\,(\varepsilon\rho)^{-(\alpha+\gamma)}}{f(\varepsilon)\,
\varepsilon^{-(\alpha+\gamma)}}
.
$$

Since $f$ is RVF of order $\alpha$, the function $f(\tau)\,\tau^{-\alpha+\gamma}$ increases for small $\tau$
and the function $f(\tau)\,\tau^{-(\alpha+\gamma)}$ decreases for small $\tau$. Therefore, we have
$$\chi_{[0,\frac{2\delta}{\varepsilon}]}(\rho)
\cdot\frac{f(\varepsilon \rho)}{f(\varepsilon)}\cdot\frac{\rho^{n-2}}{1+\rho^{q(n-\sigma)}}
\le C\bigg(\chi_{[0,1]}(\rho)\,\rho^{\alpha+n-2-\gamma}\, +\chi_{[1,+\infty[}(\rho)\,
\rho^{\alpha+n-2+\gamma-q(n-\sigma)}\bigg).
$$
Since $\alpha\ge1$, the majorant is summable at zero if $\gamma$ is 
sufficiently small. Since $\alpha\le n+1$, for small $\gamma$ the second 
exponent does not exceed
$$2n-1+\gamma-q(n-\sigma)=-1+\gamma-\frac{2n\sigma}{n-2\sigma}<-1,$$
and the majorant is summable at infinity.\hfill$\square$\medskip

\item
We recall that $F(0;z)\ge0$, and hence, for small $\delta$ and $\varepsilon$
\begin{multline*}
|y+\xi\varepsilon^{-1}\,F(\varepsilon y';\varepsilon z)\,e_m|\ge\\
\ge|y+\xi\varepsilon^{-1}\,F(0;\varepsilon z)\,e_m|-
\varepsilon^{-1}\,|F(\varepsilon y';\varepsilon z)-F(0;\varepsilon z)|\ge\\
\ge|y|-|y'|o_\delta(1)\ge\frac{|y|}{2}.
\end{multline*}
Therefore,
\begin{multline*}
\bigg|\frac{I_3 \varepsilon}{f(\varepsilon)}\bigg|\le
C\int\limits_{\R^n_+} \frac{\phi^q(y;z)}{|y|^{q(1-\sigma)+2}}\,
\widetilde\varphi^q(\varepsilon y;\varepsilon z)
\frac {F^2(\varepsilon y';\varepsilon z)}{\varepsilon\,|f(\varepsilon)|}\,dydz\\
\le C\int\limits_{\R^{n-1}}\chi^{(\varepsilon)}(|y'|;|z|)\,
\frac {F^2(\varepsilon y';\varepsilon z)}{\varepsilon\,|f(\varepsilon)|}
\int\limits_0^{\infty}\frac {\phi^{q}(y',s;z)\,ds}{(|y'|^2+s^2)^{\frac {q(1-\sigma)+2}2}}\,dy'dz.
\end{multline*}
Taking into account (\ref{estimate1}) and (\ref{estimate2}), we obtain
$$\int\limits_0^{\infty}\frac {\phi^{q}(y',s;z)\,ds}{(|y'|^2+s^2)^{\frac {q(1-\sigma)+2}2}}\le
C\int\limits_0^\infty
\frac{(|y'|^2+s^2)^{\frac{q\sigma-2}{2}}\,ds}{1+(|y'|^2+|z|^2+s^2)^{\frac{qn}2}}
\le \frac{C|y'|^{q\sigma-1}}{1+(|y'|^2+|z|^2)^{\frac{qn}{2}}},$$
and therefore,
$$\bigg|\frac{I_3 \varepsilon}{f(\varepsilon)}\bigg|\le
C\int\limits_0^{\frac{2\delta}{\varepsilon}}
\int\limits_{{\mathbb S}_\rho^{n-2}}
\frac {|y'|^{q\sigma-1}F^2(\varepsilon y';\varepsilon z)}{\varepsilon\,|f(\varepsilon)|}
\,d{\mathbb S}_\rho^{n-2}(y',z)\,\frac{d\rho}{1+\rho^{qn}}.
$$

Since $q\sigma-1>-1$, it is easy to see that $W^1_2({\mathbb S}_1^{n-2})$ is embedded into
$L_{2,w}({\mathbb S}_1^{n-2})$ with weight $w=|y'|^{q\sigma-1}$. Thus, using the Poincar\'e inequality,
we can write
\begin{multline*}
\Xint{\ \,-}\limits_{{\mathbb S}_\rho^{n-2}}\!
|y'|^{q\sigma-1}|F(y';z)|^2\,d{\mathbb S}_\rho^{n-2}(y',z)\le\\
\le C\rho^{q\sigma-1}\,\Bigg(\rho^2\cdot\!\!\Xint{\ \,-}\limits_{{\mathbb S}_\rho^{n-2}}\!
|\nabla'F(y';z)|^2\,d{\mathbb S}_\rho^{n-2}(y',z)
+\bigg(\Xint{\ \,-}\limits_{{\mathbb S}_\rho^{n-2}}\!F(y';z)\,d{\mathbb S}_\rho^{n-2}(y',z)
\bigg)^2\Bigg).
\end{multline*}
This implies, subject to (\ref{eq4}),
\begin{multline*}
\bigg|\frac{I_3 \varepsilon}{f(\varepsilon)}\bigg|\le
C\int\limits_0^{\frac{2\delta}{\varepsilon}}
\frac{(\varepsilon\rho)^2\,f_2(\varepsilon\rho)+f^2(\varepsilon\rho)}{\varepsilon\,|f(\varepsilon)|}\cdot
\frac{\rho^{n-3+q\sigma}\,d\rho}{1+\rho^{qn}}=\\
=\int\limits_0^{\frac{2\delta}{\varepsilon}}o_{\rho\varepsilon}(1)\frac{f(\varepsilon\rho)}{f(\varepsilon)}
\cdot\frac{\rho^{n-2+q\sigma}\,d\rho}{1+\rho^{qn}}=o_{\delta}(1),
\end{multline*}
and we arrive at $I_3=A_1(\varepsilon)o_{\delta}(1)$.
 \end{enumerate}

We remark also that $\frac {qn}2>n\ge\alpha-1$. This implies
$\varepsilon^{\frac{qn}{2}}=A_1(\varepsilon)o_{\varepsilon}(1)$.

Choosing $\delta>0$ sufficiently small and summing the estimates of items 1-3, we arrive at (\ref{denom}).

\section{Estimate of the numerator and derivation of (\ref{numer}) for $\alpha<n+1$}\label{numerator}

We have
$$(\phi_\varepsilon)_{x_m}=\varepsilon^{-n/2}
\phi_{y_1}(\Theta_\varepsilon(x))\,\varphi(x)+
\varepsilon^{1-n/2}\, \phi(\Theta_\varepsilon(x))\,\varphi_{x_m}(x);$$
 while for $i\ne 1$
$$(\phi_\varepsilon)_{x_i}=\varepsilon^{-n/2}
\bigg(\phi_{y_i}(\Theta_\varepsilon(x))-\phi_{y_1}(\Theta_\varepsilon(x))\,
F_{x_i}(y';z)\bigg)\varphi(x)+\varepsilon^{1-n/2}\,
\phi(\Theta_\varepsilon(x))\,\varphi_{x_i}(x).$$
Hence
\begin{multline*}
\int\limits_{\Omega}|\nabla\phi_\varepsilon(x)|^2\,dx=
\int\limits_{\Omega}
\bigg[\varepsilon^{-n}\varphi^2(x)\,|(\nabla'\phi)(\Theta_\varepsilon(x))|^2-\\
-2\,\varepsilon^{-n}\varphi^2(x)\,\phi_{y_1}(\Theta_\varepsilon(x))\,
\big\langle(\nabla'\phi)(\Theta_\varepsilon(x)),\nabla'F(y';z)\big\rangle +\vphantom{\bigg|}\\
+2\,\varepsilon^{-n+1}\phi(\Theta_\varepsilon(x))\varphi(x)
\big\langle(\nabla'\phi)(\Theta_\varepsilon(x)),\nabla'\varphi(x)\big\rangle-\vphantom{\bigg|}\\
-2\,\varepsilon^{-n+1}\phi_{y_1}(\Theta_\varepsilon(x))\,
\phi(\Theta_\varepsilon(x))\,\varphi(x)\,\big\langle\nabla'F(y';z),\nabla'\varphi(x)\big\rangle+\vphantom{\bigg|}\\
+\varepsilon^{2-n}\phi^2(\Theta_\varepsilon(x))\,|\nabla'\varphi(x)|^2
+\varepsilon^{-n}\phi^2_{y_1}(\Theta_\varepsilon(x))\,\varphi^2(x)\,|\nabla'F(y';z)|^2+\vphantom{\bigg|}\\
+\varepsilon^{-n}\varphi^2(x)\phi^2_{y_1}(\Theta_\varepsilon(x))
+\varepsilon^{2-n}\phi^2(\Theta_\varepsilon(x))\,
\varphi^2_{x_m}(x)+\vphantom{\bigg|}\\
+2\,\varepsilon^{-n+1}\phi_{y_1}(\Theta_\varepsilon(x))\,\varphi(x)\,
\phi(\Theta_\varepsilon(x))\,\varphi_{x_m}(x)\bigg]\,dx=:J_1-J_2+\dots+J_9.
\end{multline*}

\begin{enumerate}
\item $$J_1+J_7=\int\limits_{\R^n_+}
\widetilde\varphi^2(\varepsilon y;\varepsilon z)|\nabla\phi(y;z)|^2\,dydz=
\mu^2_q-\int\limits_{\R^n_+}(1-\widetilde\varphi^2(\varepsilon y;\varepsilon z))
|\nabla\phi(y;z)|^2\,dydz;$$
moreover, the second estimate in (\ref{estimate2}) gives
$$\int\limits_{\R^n_+}(1-\widetilde\varphi^2(\varepsilon y;\varepsilon z))
|\nabla\phi(y;z)|^2\,dydz\le C\int\limits_{\frac{\delta}{2\varepsilon}}^\infty 
\zeta^{-1-n}\,d\zeta=C\,\Big(\frac{\varepsilon}{\delta}\Big)^{n},$$
whence 
$$J_1+J_7=\mu^2_q+C(\delta)O(\varepsilon^{n}).$$

    \item Integrating by parts we obtain
\begin{multline*}
J_2=
\int\limits_{\R^n_+}2\,\widetilde\varphi^2(\varepsilon y;\varepsilon z)\,
\phi_{y_1}(y;z)\,\big\langle\nabla'\phi(y;z),\nabla'F(\varepsilon y;\varepsilon z)\big\rangle\,dydz=\\
=-\frac{2}{\varepsilon}\int\limits_{\R^n_+}
\bigg[\widetilde\varphi^2(\varepsilon y;\varepsilon z)\phi_{y_1}(y;z)\Delta'\phi(y;z)
+\big\langle\nabla'(\widetilde\varphi^2)(\varepsilon y;\varepsilon z),\nabla'\phi(y;z)\big\rangle\phi_{y_1}(y;z)+\\
+\widetilde\varphi^2(\varepsilon y;\varepsilon z)\big\langle\nabla'\phi(y;z),\nabla'\phi_{y_1}(y;z)\big\rangle
\bigg]\,F(\varepsilon y';\varepsilon z)\,dydz.
\end{multline*}
By (\ref{eq}), we obtain
\begin{multline*}
J_2=
\frac{2}{\varepsilon}\int\limits_{\R^n_+}\widetilde\varphi^2(\varepsilon y;\varepsilon z)
\bigg(\phi_{y_1y_1}(y;z)+
\mu^2_q\,\frac{\phi^{q-1}(y;z)}{|y|^{q(1-\sigma)}}\bigg)
\phi_{y_1}(y;z)\,F(\varepsilon y';\varepsilon z)\,dydz-\\
-\frac{2}{\varepsilon}\int\limits_{\R^n_+}
\big\langle\nabla'(\widetilde\varphi^2)(\varepsilon y;\varepsilon z),\nabla'\phi(y;z)\big\rangle
\phi_{y_1}(y;z)\,F(\varepsilon y';\varepsilon z)\,dydz-\\
-\frac{1}{\varepsilon}\int\limits_{\R^n_+}\widetilde\varphi^2(\varepsilon y;\varepsilon z)
\,(|\nabla'\phi(y;z)|^2)_{y_1}\,F(\varepsilon y';\varepsilon z)\,dydz=:
H+K_1+K_2.\end{multline*}

Now we integrate the first term by parts.
\begin{multline*}
H=\frac{1}{\varepsilon}\int\limits_{\R^n_+}\widetilde\varphi^2(\varepsilon y;\varepsilon z)
\bigg((\phi^2_{y_1}(y;z))_{y_1}+\frac{2\,\mu^2_q}{q}
\frac{(\phi^q(y;z))_{y_1}}{|y|^{q(1-\sigma)}}\bigg)
\,F(\varepsilon y';\varepsilon z)\,dydz=\\
=-\frac{1}{\varepsilon}\int\limits_{\R^{n-1}}
\widetilde\varphi^2(\varepsilon y;\varepsilon z)
\phi^2_{y_1}(0,y';z)\,F(\varepsilon y';\varepsilon z)\,dy'dz-\\
-\frac{1}{\varepsilon}\int\limits_{\R^n_+}
(\widetilde\varphi^2(\varepsilon y;\varepsilon z))_{y_1}
\,\phi^2_{y_1}(y;z)\,F(\varepsilon y';\varepsilon z)\,dydz-\\
-\frac{2\,\mu^2_q}{q\,\varepsilon}\int\limits_{\R^n_+}
(\widetilde\varphi^2(\varepsilon y;\varepsilon z))_{y_1}\frac{\phi^q(y;z)}
{|y|^{q(1-\sigma)}}\,F(\varepsilon y';\varepsilon z)\,dydz+\\
+\frac{2\,(1-\sigma)\,\mu^2_q}{\varepsilon}
\int\limits_{\R^n_+}\widetilde\varphi^2(\varepsilon y;\varepsilon z)\frac{\phi^q(y;z)}
{|y|^{q(1-\sigma)}}\frac{y_1}{|y|^2}F(\varepsilon y';\varepsilon z)\,dydz=:
-A_2(\varepsilon)+K_3+K_4+K_5.
\end{multline*}

\begin{prop}\label{prop3}. Let assumptions of Theorem \ref{Th5} hold. Then,
given $\delta$, the function $A_2(\varepsilon)$ satisfies (\ref{A2}), as $\varepsilon\to0$.
\end{prop}

\pro We claim that
\begin{multline}\lim\limits_{\varepsilon\to 0}
\frac{\varepsilon A_2(\varepsilon)}{f(\varepsilon)} = 
\omega_{m-2}\omega_{n-m-1}\times\\
\times\int\limits_0^{\infty}\rho^{n-2+\alpha}\int\limits_0^{\frac{\pi}{2}}
\cos^{m-2}(\beta)\sin^{n-m-1}(\beta)\Phi(\beta)
|\nabla \phi(0,\rho\cos(\beta);\rho\sin(\beta))|^2\, d\beta d\rho.
\label{eq7}\end{multline}

To prove this we apply the Lebesgue theorem. We have, similarly to Proposition \ref{prop2}, 
\begin{equation}\label{eq8}
\frac{\varepsilon A_2(\varepsilon)}{\omega_{m-2}\omega_{n-m-1}f(\varepsilon)}=
\int\limits_0^{\infty}\int\limits_0^{\infty}
\widetilde{\varphi}^2(\varepsilon r;\varepsilon t)\,r^{m-2}t^{n-m-1}\,
\frac{f_1(\varepsilon r;\varepsilon t)}{f(\varepsilon)}\,|\nabla \phi(0,r;t)|^2\,drdt.
\end{equation}

Passing to the polar coordinates, we see that the integrand converges to 
that in (\ref{eq7}) for all $\rho$ and $\beta$. Now we produce a summable majorant.
By (\ref{estimate1}) and (\ref{estimate2}), we have
$$|\nabla \phi(0,\rho\cos(\beta);\rho\sin(\beta))|\le \frac{C}{1+\rho^n}.
$$
Therefore for $\gamma>0$, similarly to Proposition \ref{prop2}, we estimate the integrand by
$$C\chi_{[0,\frac{2\delta}{\varepsilon}]}(\rho)\,\frac{f(\varepsilon\, \rho)}{f(\varepsilon)}\,
\frac{\rho^{n-2}}{1+\rho^{2n}}\le\\
C\Bigl(\chi_{[0,1]}(\rho)\,\rho^{\alpha+n-2-\gamma}+
\chi_{[1,\infty[}(\rho)\,\rho^{\alpha-n-2+\gamma}\Bigr).
$$
Since $\alpha<n+1$, this provides a summable majorant for sufficiently small $\gamma$.
\hfill$\square$\medskip

Now we estimate all remaining terms in $J_2$. Since the functions
$\widetilde\varphi$ and $\phi$ are radially symmetric w.r.t $y'$ and $z$,
integrating by parts in $K_2$ we have
\begin{multline*}
K_1+K_2+K_3=
-\frac{1}{\varepsilon}\int\limits_{\R^n_+}
F(\varepsilon y';\varepsilon z)\,
\bigg[2\,\big\langle\nabla'(\widetilde\varphi^2)(\varepsilon y;\varepsilon z),\nabla'\phi(y;z)\big\rangle
\phi_{y_1}(y;z)-\\
-(\widetilde\varphi^2(\varepsilon y;\varepsilon z))_{y_1}\,|\nabla'\phi(y;z)|^2+
(\widetilde\varphi^2(\varepsilon y;\varepsilon z))_{y_1}\,\phi^2_{y_1}(y;z)\bigg]\,dydz=
\end{multline*}
\begin{multline*}
=-\frac 1{\varepsilon}
\int\limits_0^{\infty}\int\limits_0^{\infty}\int\limits_0^{\infty}
r^{m-2}t^{n-m-1}f_1(\varepsilon r;\varepsilon t)\times\\
\times\Big[2\big((\widetilde\varphi^2(\varepsilon y_1,\varepsilon r;\varepsilon t))_r\,\phi_r+
(\widetilde\varphi^2(\varepsilon y_1,\varepsilon r;\varepsilon t))_t\,\phi_t\big)\phi_{y_1}-\\
-(\widetilde\varphi^2(\varepsilon y_1,\varepsilon r;\varepsilon t))_{y_1}\,\big(\phi_r^2+\phi_t^2\big)
+(\widetilde\varphi^2(\varepsilon y_1,\varepsilon r;\varepsilon t))_{y_1}\,\phi_{y_1}^2\Big]
\,dy_1\,drdt.
\end{multline*}
Using the assumption (\ref{eq11}) and the estimate (\ref{estimate2}), we obtain
\begin{multline*}
|K_1+K_2+K_3|\le\frac{C}{\delta}\iint\limits_{\frac{\delta}{2\varepsilon}\le
\sqrt{\rho^2+y_1^2}\le\frac{2\delta}{\varepsilon}}
\!\!\!|\nabla\phi|^2\,\rho^{n-2}|f(\varepsilon\rho)|\,dy_1\,d\rho\le\\
\le\frac{C}{\delta}\iint\limits_{\frac{\delta}{2\varepsilon}\le 
\sqrt{\rho^2+s^2}\le\frac{2\delta}{\varepsilon}}^{\vphantom{1^1}}
\frac {\rho^{n-2}\,|f(\varepsilon\rho)|}{{(\rho^2+s^2)}^n}\,ds\,d\rho=
\frac {C\varepsilon^n}{\delta}\iint\limits_{\frac{\delta}{2}\le 
\sqrt{\mathfrak r^2+\mathfrak s^2}\le 2\delta}
\frac {\mathfrak r^{n-2}\,|f(\mathfrak r)|}{(\mathfrak r^2+\mathfrak s^2)^n}\,d\mathfrak s\,d\mathfrak r
=C(\delta)\cdot\varepsilon^n.
\end{multline*}

In a similar way,
$$K_4=-\frac{2\,\mu^2_q}{q\,\varepsilon}
\int\limits_0^{\infty}\int\limits_0^{\infty}\int\limits_0^{\infty}
r^{m-2}t^{n-m-1}f_1(\varepsilon r;\varepsilon t)\cdot
(\widetilde\varphi^2(\varepsilon y_1,\varepsilon r;\varepsilon t))_{y_1}\,
\frac{\phi^q(y_1,r;t)}{(r^2+y_1^2)^{\frac{q(1-\sigma)}2}}\,dy_1\,drdt,
$$
and therefore,
\begin{multline*}
|K_4|\le\frac{C}{\delta}\iint\limits_{\frac{\delta}{2\varepsilon}\le
\sqrt{\rho^2+y_1^2}\le\frac{2\delta}{\varepsilon}}
\frac{y_1^{q\sigma}\rho^{n-2}|f(\varepsilon\rho)|}{(\rho^2+y_1^2)^{\frac{qn}{2}}}\,dy_1\,d\rho=\\
=\frac {C\varepsilon^{q(n-\sigma)-n}}{\delta}\!\!\iint\limits_{\frac{\delta}{2}\le 
\sqrt{\mathfrak r^2+\mathfrak s^2}\le 2\delta}^{\vphantom{1^1}}
\frac {\mathfrak s^{q\sigma} \mathfrak r^{n-2}\,|f(\mathfrak r)|}
{(\mathfrak r^2+\mathfrak s^2)^{\frac{qn}{2}}}\,d\mathfrak s\,d\mathfrak r
=C(\delta)\cdot\varepsilon^{q(n-\sigma)-n}=o(\varepsilon^n)
\end{multline*}
(the last relation follows from
$q(n-\sigma)-2n=\frac {2n\sigma}{n-2\sigma}>0$).

Finally, the integral $K_5$ can be estimated in the same way as $I_2$ in Section \ref{denominator}. 
This gives, as $\varepsilon\to0$,
$$K_5\sim \frac {2\mu^2_q}{q}A_1(\varepsilon).$$
Thus, 
$$J_2=-A_2(\varepsilon)+C(\delta)O(\varepsilon^n)+
\frac {2\mu^2_q}{q}A_1(\varepsilon)(1+o_\varepsilon(1)).$$

    \item 
By the estimate (\ref{estimate2}), we obtain
\begin{multline*}
|J_3+J_9|\le 2\,\varepsilon\int\limits_{\R^n_+}\phi(y;z)\,
|\nabla \phi(y;z)|\,\widetilde\varphi(\varepsilon y;\varepsilon z)\,
|\nabla\widetilde\varphi(\varepsilon y;\varepsilon z)|\,dydz\le\\
\le\frac {C\,\varepsilon}{\delta}
\int\limits_{\frac{\delta}{2\varepsilon}}^{\frac{2\delta}{\varepsilon}}
\zeta^{-n}\,d\zeta=C\,\Big(\frac{\varepsilon}{\delta}\Big)^{n}.
\end{multline*}

    \item 
Using the previous estimate, we obviously get
\begin{multline*}
|J_4|\le 2\,\varepsilon\int\limits_{\R^n_+} \phi(y;z)\,|\nabla\phi(y;z)|\,
\widetilde\varphi(\varepsilon y;\varepsilon z)\,|\nabla\widetilde\varphi(\varepsilon y;\varepsilon z)|\,
|\nabla'F(\varepsilon y';\varepsilon z)|\,dydz\le\\
\le C\,\varepsilon\int\limits_{\R^n_+} \phi(y;z)\,|\nabla\phi(y;z)|\,
\widetilde\varphi(\varepsilon y;\varepsilon z)\,|\nabla\widetilde\varphi(\varepsilon y;\varepsilon z)|
\,dydz\le C\,\Big(\frac{\varepsilon}{\delta}\Big)^{n}.
\end{multline*}

\item 
In a similar way,
$$|J_5+J_8|=\varepsilon^{2}\int\limits_{\R^n_+}\phi^2(y;z)
\,|\nabla\widetilde\varphi(\varepsilon y;\varepsilon z)|^2\,dydz
\le \frac {C\,\varepsilon^2}{\delta^2}
\int\limits_{\frac{\delta}{2\varepsilon}}^{2\frac{\delta}{\varepsilon}}
\zeta^{1-n}\,d\zeta=C\,\Big(\frac{\varepsilon}{\delta}\Big)^{n}.$$

\item Finally, relations (\ref{estimate2}) and (\ref{eq4}) imply
\begin{multline*}
J_6=\int\limits_{\R^n_+}
\phi^2_{y_1}(y;z)\,\widetilde\varphi^2(\varepsilon y;\varepsilon z)\,
|\nabla'F(\varepsilon y';\varepsilon z)|^2\,dydz\le \\
\le C\int\limits_{0}^
{\frac{2\delta}{\varepsilon}}\rho^{n-2}\,f_2(\varepsilon\rho)
\int\limits_0^\infty\frac {dy_1}{(1+\rho^2+y_1^2)^n}\,d\rho=\\
\le C\int\limits_{0}^{\frac{2\delta}{\varepsilon}}\frac 
{\rho^{n-2}\,f_2(\varepsilon\rho)}{(1+\rho^2)^{n-1/2}}\,d\rho
=\frac {o_\delta(1)}{\varepsilon}\int\limits_{0}^{\frac{2\delta}{\varepsilon}}
\frac {\rho^{n-3}\,|f(\varepsilon\rho)|}{(1+\rho^2)^{n-1/2}}\,d\rho.
\end{multline*}

The last integral can be estimated in the same way as in Proposition \ref{prop3}.
This gives
$$J_5=o_\delta(1)A_2(\varepsilon).$$
\end{enumerate}

We remark also that $\varepsilon^n=A_2(\varepsilon)o_{\varepsilon}(1)$.

Choosing $\delta>0$ sufficiently small and summing the estimates of items 1-6
we obtain (\ref{numer}).

\section{Derivation of (\ref{numer}) for $\alpha=n+1$}\label{limit case}

We underline that the assumption $\alpha<n+1$ was used in the previous section only
in the proof of Proposition \ref{prop3}. Also the assumption (\ref{eq2}) was used only
to ensure the positivity of the integral in (\ref{eq7}). So, we need only to prove
the following fact.

\begin{prop}\label{prop4} Let assumptions of Theorem \ref{Th6} hold. Then,
given $\delta$, the function $A_2(\varepsilon)$ satisfies (\ref{A2}), as $\varepsilon\to0$.
\end{prop}

\pro By (\ref{estimate2}), there exists $M>0$ such that
\begin{equation}
|\nabla \phi(0,r;t)|=\frac {M+o_\rho(1)}{\rho^n},\qquad \rho=\sqrt{r^2+t^2}\to\infty.
\label{qqq}\end{equation}
We split the integral (\ref{eq8}) into three parts:
\begin{multline*}
\frac{A_2(\varepsilon)}{\omega_{m-2}\omega_{n-m-1}\varepsilon^n}=
\bigg(\,\iint\limits_{\sqrt{r^2+t^2}\le R}+\iint\limits_{R\le\sqrt{r^2+t^2}\le\frac{\delta}{2\varepsilon}}+
\iint\limits_{\frac{\delta}{2\varepsilon}\le\sqrt{r^2+t^2}\le\frac{2\delta}{\varepsilon}}\bigg)
\widetilde{\varphi}^2(\varepsilon r;\varepsilon t)\times\\
\times r^{m-2}t^{n-m-1}\,\frac{f_1(\varepsilon r;\varepsilon t)}{\varepsilon^{n+1}}
\,|\nabla \phi(0,r;t)|^2\,drdt=:L_1+L_2+L_3.\vphantom{\Bigg|}
\end{multline*}
The relation (\ref{qqq}) implies that, as $R\to\infty$,
\begin{multline*}
L_2=(M+o_R(1))^2\times\\
\times\int\limits_R^{\frac{\delta}{2\varepsilon}}
\int\limits_0^{\frac{\pi}{2}}\cos^{m-2}(\beta)\sin^{n-m-1}(\beta)\,
\frac{f_1(\varepsilon \rho\cos(\beta);\varepsilon \rho\sin(\beta))}
{\varepsilon^{n+1}\rho^{n+2}}\,d\beta\,d\rho=\\
=\frac{(M^2+o_R(1))\,\omega_{n-2}}{\omega_{m-2}\omega_{n-m-1}}
\int\limits_R^{\frac{\delta}{2\varepsilon}}\frac{f(\varepsilon \rho)}{\varepsilon^{n+1}\rho^{n+2}}\,d\rho=
\frac{(M^2+o_R(1))\,\omega_{n-2}}{\omega_{m-2}\omega_{n-m-1}}
\int\limits_{R\varepsilon}^{\delta/2}\frac {f(\mathfrak r)}{\mathfrak r^{n+2}}\,d\mathfrak r.
\end{multline*}

Further, the assumption (\ref{eq11}) implies
$$\left|\frac {\varepsilon^{n+1}L_1}{f(\varepsilon)}\right|\le C\int\limits_0^R
\frac {f(\varepsilon\rho)}{f(\varepsilon)}\,\rho^{n-2}\,d\rho.$$
For given $R$ we can pass to the limit under the integral sign. This provides
$L_1=O\big(\frac{f(\varepsilon)}{\varepsilon^{n+1}}\big)$, as $\varepsilon\to0$.

On the another hand, divergence of the integral 
$\int_0^\delta\frac {f(\mathfrak r)}{\mathfrak r^{n+2}}\,d\mathfrak r$ implies that
for
arbitrary large $N$ we have, as $\varepsilon$ is sufficiently small,
\begin{equation}
\frac {\varepsilon^{n+1}}{f(\varepsilon)}\,\int\limits_{R\varepsilon}^{\delta/2}
\frac {f(\mathfrak r)}{\mathfrak r^{n+2}}\,d\mathfrak r\ge
\frac {\varepsilon^{n+1}}{f(\varepsilon)}\,\int\limits_{R\varepsilon}^{N\varepsilon}
\frac {f(\mathfrak r)}{\mathfrak r^{n+2}}\,d\mathfrak r
=\int\limits_R^N\frac{f(\varepsilon\rho)}
{f(\varepsilon)\rho^{n+1}}\,\frac{d\rho}{\rho}=\ln(N/R)\cdot(1+o_\varepsilon(1)),
\label{svf}\end{equation}
and thus, $L_1=o(L_2)$, as $\varepsilon\to0$.

Finally, as $\varepsilon\to0$,
$$|L_3|\le C\int\limits_{\frac{\delta}{2\varepsilon}}^{\frac{2\delta}{\varepsilon}}
\frac{f(\varepsilon \rho)}{\varepsilon^{n+1}\rho^{n+2}}\,d\rho
=C\int\limits_{\delta/2}^{2\delta}\frac {f(\mathfrak r)}{\mathfrak r^{n+2}}\,d\mathfrak r
=C(\delta)=o(L_2).$$

It remains to note that for given $R$ and $\delta$ 
$$\int\limits_{R\varepsilon}^{\delta/2}\frac {f(\mathfrak r)}{\mathfrak r^{n+2}}\,d\mathfrak r=
\int\limits_{\varepsilon}^{\delta}\frac {f(\mathfrak r)}{\mathfrak r^{n+2}}\,d\mathfrak r+O(1)\sim
\int\limits_{\varepsilon}^{\delta}\frac {f(\mathfrak r)}{\mathfrak r^{n+2}}\,d\mathfrak r, 
\qquad \varepsilon\to0,$$
and we arrive at
$$A_2(\varepsilon)\sim M^2\omega_{n-2}\varepsilon^{n}
\int\limits_{\varepsilon}^{\delta}\frac {f(\mathfrak r)}{\mathfrak r^{n+2}}\,d\mathfrak r.$$
\hfill$\square$

\end{document}